\documentclass[10pt]{article}
\usepackage[utf8]{inputenc}
\usepackage[T1]{fontenc}
\usepackage{lmodern}
\usepackage{amsmath,amssymb,amsthm}
\usepackage{array,booktabs,graphicx}
\usepackage[margin=1.5in]{geometry}
\usepackage{indentfirst}
\usepackage{tikz}
\usetikzlibrary{arrows,arrows.meta,shapes.geometric,positioning,calc}
\usepackage[font=small,labelsep=period]{caption}
\usepackage{subcaption}
\usepackage[export]{adjustbox}
\usepackage{placeins}
\usepackage[hidelinks]{hyperref}
\newcommand{\Z}{\mathbb Z}

\newcommand{\E}{\mathbb E}
\DeclareMathOperator{\Aut}{Aut}
\DeclareMathOperator{\lca}{lca}
\newtheorem{theorem}{Theorem}
\newtheorem{lemma}{Lemma}
\newtheorem{proposition}{Proposition}
\newtheorem{corollary}{Corollary}
\theoremstyle{definition}
\newtheorem{remark}{Remark}

\newcommand{\correspondingauthor}{\thanks{Corresponding author}}

\newcommand{\GallaiFigureOne}{%
\begin{figure}[htbp]
\centering
\begin{adjustbox}{max width=\linewidth}
\begin{tikzpicture}[scale=0.85, baseline=0]
    \tikzset{
        vnode/.style={circle, fill=black, inner sep=1.2pt},
        c1/.style={red, thick},
        c2/.style={blue, thick},
        c3/.style={green!70!black, thick}
    }

    \node[vnode, label={135:$v_2$}] (w0) at (-4, 1.5) {};
    \node[vnode, label={225:$v_1$}] (w1) at (-4, -1.5) {};
    \node[vnode, label={315:$v_4$}] (w2) at (-1, -1.5) {};
    \node[vnode, label={45:$v_3$}] (w3) at (-1, 1.5) {};

    \draw[c1] (w1) -- (w0); 
    \draw[c1] (w0) -- (w3); 
    \draw[c1] (w3) -- (w2); 
    
    \draw[c2] (w1) -- (w2); 
    \draw[c2] (w2) -- (w0); 
    \draw[c2] (w1) -- (w3); 

    \node[font=\large] at (-2.5, -2.5) {$\Pi$};

    \node[vnode, label={90:$v_2$}] (v0) at (2.5, 1.5) {};
    \node[vnode, label={225:$v_1$}] (v1) at (1.0, -1.5) {};
    \node[vnode, label={315:$v_3$}] (v2) at (4.0, -1.5) {};

    \draw[c1] (v1) -- (v0); 
    \draw[c2] (v0) -- (v2); 
    \draw[c3] (v2) -- (v1); 
    
    \node[font=\large] at (2.5, -2.5) {$\Delta$};
    
\end{tikzpicture}
\end{adjustbox}
\caption{$2$-graph $\Pi$ and $3$-graph $\Delta$}
\label{f1}
\end{figure}
}

\newcommand{\GallaiFigureTwo}{%
\begin{figure}[htbp]
\centering
\begin{subfigure}[c]{0.32\textwidth}
\centering
\begin{adjustbox}{max width=\linewidth}
\begin{tikzpicture}[scale=0.85, baseline=0]
    \tikzset{
        tnode/.style={circle, draw=black, fill=white, inner sep=0.5pt, font=\scriptsize\sffamily, minimum size=13pt},
        leaf/.style={circle, fill=black, inner sep=1.2pt}
    }

    \node[circle, draw=black, fill=white, inner sep=0pt, minimum size=3pt,
        label={[font=\scriptsize]left:$s$}] (vs) at (-0.7,0) {};
    
    \node[tnode] (v0) at (0,0) {1};
    \draw[densely dashed] (vs) -- (v0);

    \node[tnode] (N1_0) at (90:1.1) {2};
    \node[tnode] (N1_1) at (270:1.1) {2};
    \draw (v0) -- (N1_0);
    \draw (v0) -- (N1_1);

    \node[leaf, label={45:$\omega_0$}] (L0) at (45:2.0) {};
    \node[leaf, label={135:$\omega_1$}] (L1) at (135:2.0) {};
    \node[leaf, label={225:$\omega_2$}] (L2) at (225:2.0) {};
    \node[leaf, label={315:$\omega_3$}] (L3) at (315:2.0) {};

    \draw (N1_0) -- (L0);
    \draw (N1_0) -- (L1);
    \draw (N1_1) -- (L2);
    \draw (N1_1) -- (L3);
\end{tikzpicture}
\end{adjustbox}
\caption{Game tree $T$ ($m=2$)}
\label{fig:tree_m2}
\end{subfigure}
\hfill
\begin{subfigure}[c]{0.32\textwidth}
\centering
\begin{adjustbox}{max width=\linewidth}
\begin{tikzpicture}[scale=0.85, baseline=0]
    \tikzset{
        c1/.style={red, thick},
        c2/.style={blue, thick},
        vnode/.style={circle, fill=black, inner sep=1.2pt}
    }

    \node[vnode, label={45:$\omega_0$}] (V0) at (45:2.0) {};
    \node[vnode, label={135:$\omega_1$}] (V1) at (135:2.0) {};
    \node[vnode, label={225:$\omega_2$}] (V2) at (225:2.0) {};
    \node[vnode, label={315:$\omega_3$}] (V3) at (315:2.0) {};

    \draw[c2] (V0) -- (V1);
    \draw[c2] (V2) -- (V3);

    \draw[c1] (V0) -- (V2);
    \draw[c1] (V0) -- (V3);
    \draw[c1] (V1) -- (V2);
    \draw[c1] (V1) -- (V3);
\end{tikzpicture}
\end{adjustbox}
\caption{2-colored $K_4$}
\label{fig:k4_colored}
\end{subfigure}
\hfill
\begin{subfigure}[c]{0.32\textwidth}
\centering
\begin{adjustbox}{max width=\linewidth}
\begin{tikzpicture}[scale=0.85, baseline=0]
    \tikzset{
        c1/.style={red, thick},
        c2/.style={blue, thick},
        vnode/.style={circle, fill=black, inner sep=1.2pt}
    }

    \node[vnode, label={45:$\omega_0$}] (V0) at (45:2.0) {};
    \node[vnode, label={135:$\omega_2$}] (V2) at (135:2.0) {};
    \node[vnode, label={225:$\omega_1$}] (V1) at (225:2.0) {};
    \node[vnode, label={315:$\omega_3$}] (V3) at (315:2.0) {};

    \draw[c1] (V0) -- (V2);
    \draw[c1] (V2) -- (V1);
    \draw[c1] (V1) -- (V3);
    \draw[c1] (V3) -- (V0);

    \draw[c2] (V0) -- (V1);
    \draw[c2] (V2) -- (V3);
\end{tikzpicture}
\end{adjustbox}
\caption{2-circulant representation}
\label{fig:circulant_m2}
\end{subfigure}

\caption{The binary tree $T$ of depth $m=2$ with branching numbers $k_1=k_2=2$, the corresponding 2-graph, and its circular representation $(1;2)_4$, obtained by ordering the leaves as $(\omega_0,\omega_2,\omega_1,\omega_3)$. The auxiliary vertex $s$ is not a game position; play starts at $v_0$, controlled by player 1.}
\label{f2}
\end{figure}
}

\newcommand{\GallaiFigureThree}{%
\begin{figure}[htbp]
\centering
\begin{subfigure}[c]{0.32\textwidth}
\centering
\begin{adjustbox}{max width=\linewidth}
\begin{tikzpicture}[scale=0.85, baseline=0]
    \tikzset{
        tnode/.style={circle, draw=black, fill=white, inner sep=0.5pt, font=\scriptsize\sffamily, minimum size=13pt},
        leaf/.style={circle, fill=black, inner sep=1.2pt}
    }

    \node[circle, draw=black, fill=white, inner sep=0pt, minimum size=3pt,
        label={[font=\scriptsize]left:$s$}] (vs) at (-0.7, 0) {};
    
    \node[tnode] (v0) at (0,0) {1};
    \draw[densely dashed] (vs) -- (v0);

    \node[tnode] (N1_0) at (0, 1.2) {2};
    \node[tnode] (N1_1) at (0, -1.2) {2};
    \draw (v0) -- (N1_0);
    \draw (v0) -- (N1_1);

    \node[leaf, label={0:$\omega_0$}] (L0) at (0.8, 2.2) {};
    \node[leaf, label={90:$\omega_1$}] (L1) at (0, 2.2) {};
    \node[leaf, label={180:$\omega_2$}] (L2) at (-0.8, 2.2) {};
    
    \node[leaf, label={180:$\omega_3$}] (L3) at (-0.8, -2.2) {};
    \node[leaf, label={270:$\omega_4$}] (L4) at (0, -2.2) {};
    \node[leaf, label={0:$\omega_5$}] (L5) at (0.8, -2.2) {};

    \draw (N1_0) -- (L0);
    \draw (N1_0) -- (L1);
    \draw (N1_0) -- (L2);
    
    \draw (N1_1) -- (L3);
    \draw (N1_1) -- (L4);
    \draw (N1_1) -- (L5);
\end{tikzpicture}
\end{adjustbox}
\caption{Game tree $T$}
\label{fig:tree_m3}
\end{subfigure}
\hfill
\begin{subfigure}[c]{0.32\textwidth}
\centering
\begin{adjustbox}{max width=\linewidth}
\begin{tikzpicture}[scale=0.85, baseline=0]
    \tikzset{
        c1/.style={red, thick},
        c2/.style={blue, thick},
        vnode/.style={circle, fill=black, inner sep=1.2pt}
    }

    \node[vnode, label={30:$\omega_0$}] (V0) at (30:2.2) {};
    \node[vnode, label={90:$\omega_1$}] (V1) at (90:2.2) {};
    \node[vnode, label={150:$\omega_2$}] (V2) at (150:2.2) {};
    \node[vnode, label={210:$\omega_3$}] (V3) at (210:2.2) {};
    \node[vnode, label={270:$\omega_4$}] (V4) at (270:2.2) {};
    \node[vnode, label={330:$\omega_5$}] (V5) at (330:2.2) {};

    \draw[c2] (V0) -- (V1); \draw[c2] (V1) -- (V2); \draw[c2] (V2) -- (V0);
    \draw[c2] (V3) -- (V4); \draw[c2] (V4) -- (V5); \draw[c2] (V5) -- (V3);

    \draw[c1] (V0) -- (V3); \draw[c1] (V0) -- (V4); \draw[c1] (V0) -- (V5);
    \draw[c1] (V1) -- (V3); \draw[c1] (V1) -- (V4); \draw[c1] (V1) -- (V5);
    \draw[c1] (V2) -- (V3); \draw[c1] (V2) -- (V4); \draw[c1] (V2) -- (V5);
\end{tikzpicture}
\end{adjustbox}
\caption{2-colored $K_6$}
\label{fig:k6_colored}
\end{subfigure}
\hfill
\begin{subfigure}[c]{0.32\textwidth}
\centering
\begin{adjustbox}{max width=\linewidth}
\begin{tikzpicture}[scale=0.85, baseline=0]
    \tikzset{
        c1/.style={red, thick},
        c2/.style={blue, thick},
        vnode/.style={circle, fill=black, inner sep=1.2pt}
    }

    \node[vnode, label={30:$\omega_0$}] (V0) at (30:2.2) {};
    \node[vnode, label={90:$\omega_3$}] (V3) at (90:2.2) {};
    \node[vnode, label={150:$\omega_1$}] (V1) at (150:2.2) {};
    \node[vnode, label={210:$\omega_4$}] (V4) at (210:2.2) {};
    \node[vnode, label={270:$\omega_2$}] (V2) at (270:2.2) {};
    \node[vnode, label={330:$\omega_5$}] (V5) at (330:2.2) {};

    \draw[c1] (V0) -- (V3);
    \draw[c1] (V3) -- (V1);
    \draw[c1] (V1) -- (V4);
    \draw[c1] (V4) -- (V2);
    \draw[c1] (V2) -- (V5);
    \draw[c1] (V5) -- (V0);

    \draw[c1] (V0) -- (V4);
    \draw[c1] (V3) -- (V2);
    \draw[c1] (V1) -- (V5);

    \draw[c2] (V0) -- (V1);
    \draw[c2] (V1) -- (V2);
    \draw[c2] (V2) -- (V0);
    
    \draw[c2] (V3) -- (V4);
    \draw[c2] (V4) -- (V5);
    \draw[c2] (V5) -- (V3);
\end{tikzpicture}
\end{adjustbox}
\caption{2-circulant representation}
\label{fig:circulant_m3}
\end{subfigure}

\caption{Tree $T$ of depth $m=2$ with branching numbers $k_1=2$ and $k_2=3$, the corresponding 2-graph, and its circular representation $(1,3;2)_6$, obtained by ordering the leaves as $(\omega_0,\omega_3,\omega_1,\omega_4,\omega_2,\omega_5)$.}
\label{f3}
\end{figure}
}

\newcommand{\GallaiFigureFour}{%
\begin{figure}[htbp]
\centering
\begin{subfigure}[c]{0.32\textwidth}
\centering
\begin{adjustbox}{max width=\linewidth}
\begin{tikzpicture}[scale=0.8, baseline=0]
    \tikzset{
        tnode/.style={circle, draw=black, fill=white, inner sep=0.5pt, font=\scriptsize\sffamily, minimum size=13pt},
        leaf/.style={circle, fill=black, inner sep=1.2pt}
    }

    \node[circle, draw=black, fill=white, inner sep=0pt, minimum size=3pt,
        label={[font=\scriptsize]left:$s$}] (vs) at (-1.0, 0) {};
    
    \node[tnode] (v0) at (0,0) {1};
    \draw[densely dashed] (vs) -- (v0);

    \node[tnode] (N1_0) at (90:1.2) {2};
    \node[tnode] (N1_1) at (210:1.2) {2};
    \node[tnode] (N1_2) at (330:1.2) {2};
    \draw (v0) -- (N1_0);
    \draw (v0) -- (N1_1);
    \draw (v0) -- (N1_2);

    \node[leaf, label={45:$\omega_0$}] (L0) at (60:2.2) {};
    \node[leaf, label={135:$\omega_1$}] (L1) at (120:2.2) {};
    
    \node[leaf, label={180:$\omega_2$}] (L2) at (180:2.2) {};
    \node[leaf, label={270:$\omega_3$}] (L3) at (240:2.2) {};
    
    \node[leaf, label={270:$\omega_4$}] (L4) at (300:2.2) {};
    \node[leaf, label={0:$\omega_5$}] (L5) at (0:2.2) {};

    \draw (N1_0) -- (L0);
    \draw (N1_0) -- (L1);
    
    \draw (N1_1) -- (L2);
    \draw (N1_1) -- (L3);
    
    \draw (N1_2) -- (L4);
    \draw (N1_2) -- (L5);
\end{tikzpicture}
\end{adjustbox}
\caption{Game tree $T$}
\label{fig:tree_m4}
\end{subfigure}
\hfill
\begin{subfigure}[c]{0.32\textwidth}
\centering
\begin{adjustbox}{max width=\linewidth}
\begin{tikzpicture}[scale=0.8, baseline=0]
    \tikzset{
        c1/.style={red, thick},
        c2/.style={blue, thick},
        vnode/.style={circle, fill=black, inner sep=1.2pt}
    }

    \node[vnode, label={60:$\omega_0$}] (V0) at (60:2.2) {};
    \node[vnode, label={120:$\omega_1$}] (V1) at (120:2.2) {};
    \node[vnode, label={180:$\omega_2$}] (V2) at (180:2.2) {};
    \node[vnode, label={240:$\omega_3$}] (V3) at (240:2.2) {};
    \node[vnode, label={300:$\omega_4$}] (V4) at (300:2.2) {};
    \node[vnode, label={0:$\omega_5$}] (V5) at (0:2.2) {};

    \draw[c2] (V0) -- (V1); 
    \draw[c2] (V2) -- (V3); 
    \draw[c2] (V4) -- (V5);

    \draw[c1] (V0) -- (V2); \draw[c1] (V0) -- (V3); \draw[c1] (V0) -- (V4); \draw[c1] (V0) -- (V5);
    \draw[c1] (V1) -- (V2); \draw[c1] (V1) -- (V3); \draw[c1] (V1) -- (V4); \draw[c1] (V1) -- (V5);
    \draw[c1] (V2) -- (V4); \draw[c1] (V2) -- (V5);
    \draw[c1] (V3) -- (V4); \draw[c1] (V3) -- (V5);
\end{tikzpicture}
\end{adjustbox}
\caption{2-colored $K_6$}
\label{fig:k6_colored_4}
\end{subfigure}
\hfill
\begin{subfigure}[c]{0.32\textwidth}
\centering
\begin{adjustbox}{max width=\linewidth}
\begin{tikzpicture}[scale=0.8, baseline=0]
    \tikzset{
        c1/.style={red, thick},
        c2/.style={blue, thick},
        vnode/.style={circle, fill=black, inner sep=1.2pt}
    }

    \node[vnode, label={60:$\omega_0$}] (V0) at (60:2.2) {};
    \node[vnode, label={120:$\omega_2$}] (V2) at (120:2.2) {};
    \node[vnode, label={180:$\omega_4$}] (V4) at (180:2.2) {};
    \node[vnode, label={240:$\omega_1$}] (V1) at (240:2.2) {};
    \node[vnode, label={300:$\omega_3$}] (V3) at (300:2.2) {};
    \node[vnode, label={0:$\omega_5$}] (V5) at (0:2.2) {};

    \draw[c1] (V0) -- (V2);
    \draw[c1] (V2) -- (V4);
    \draw[c1] (V4) -- (V1);
    \draw[c1] (V1) -- (V3);
    \draw[c1] (V3) -- (V5);
    \draw[c1] (V5) -- (V0);

    \draw[c1] (V0) -- (V4);
    \draw[c1] (V4) -- (V3);
    \draw[c1] (V3) -- (V0);
    
    \draw[c1] (V2) -- (V1);
    \draw[c1] (V1) -- (V5);
    \draw[c1] (V5) -- (V2);

    \draw[c2] (V0) -- (V1);
    \draw[c2] (V2) -- (V3);
    \draw[c2] (V4) -- (V5);
\end{tikzpicture}
\end{adjustbox}
\caption{2-circulant representation}
\label{fig:circulant_m4}
\end{subfigure}

\caption{Tree $T$ of depth $m=2$ with branching numbers $k_1=3$ and $k_2=2$, the corresponding 2-graph, and its circular representation $(1,2;3)_6$, obtained by ordering the leaves as $(\omega_0,\omega_2,\omega_4,\omega_1,\omega_3,\omega_5)$.}
\label{f4}
\end{figure}
}

\newcommand{\GallaiFigureFive}{%
\begin{figure}[htbp]
\centering
\begin{subfigure}[c]{0.32\textwidth}
\centering
\begin{adjustbox}{max width=\linewidth}
\begin{tikzpicture}[scale=0.75, baseline=0]
    \tikzset{
        tnode/.style={circle, draw=black, fill=white, inner sep=0.5pt, font=\scriptsize\sffamily, minimum size=13pt},
        leaf/.style={circle, fill=black, inner sep=1.2pt}
    }

    \node[circle, draw=black, fill=white, inner sep=0pt, minimum size=3pt,
        label={[font=\scriptsize]left:$s$}] (vs) at (-0.7, 0) {};
    
    \node[tnode] (v0) at (0,0) {1};
    \draw[densely dashed] (vs) -- (v0);

    \node[tnode] (N1_0) at (0, 1.2) {2};
    \node[tnode] (N1_1) at (0, -1.2) {2};
    \draw (v0) -- (N1_0);
    \draw (v0) -- (N1_1);

    \node[tnode] (N2_0) at (-0.8, 2.3) {3};
    \node[tnode] (N2_1) at (0.8, 2.3) {3};
    \node[tnode] (N2_2) at (-0.8, -2.3) {3};
    \node[tnode] (N2_3) at (0.8, -2.3) {3};

    \draw (N1_0) -- (N2_0); \draw (N1_0) -- (N2_1);
    \draw (N1_1) -- (N2_2); \draw (N1_1) -- (N2_3);

    \node[leaf, label={90:$\omega_0$}] (L0) at (-1.2, 3.2) {};
    \node[leaf, label={90:$\omega_1$}] (L1) at (-0.4, 3.2) {};
    \node[leaf, label={90:$\omega_2$}] (L2) at (0.4, 3.2) {};
    \node[leaf, label={90:$\omega_3$}] (L3) at (1.2, 3.2) {};

    \node[leaf, label={270:$\omega_4$}] (L4) at (-1.2, -3.2) {};
    \node[leaf, label={270:$\omega_5$}] (L5) at (-0.4, -3.2) {};
    \node[leaf, label={270:$\omega_6$}] (L6) at (0.4, -3.2) {};
    \node[leaf, label={270:$\omega_7$}] (L7) at (1.2, -3.2) {};

    \draw (N2_0) -- (L0); \draw (N2_0) -- (L1);
    \draw (N2_1) -- (L2); \draw (N2_1) -- (L3);
    \draw (N2_2) -- (L4); \draw (N2_2) -- (L5);
    \draw (N2_3) -- (L6); \draw (N2_3) -- (L7);
\end{tikzpicture}
\end{adjustbox}
\caption{Game tree $T$ ($m=3$)}
\label{fig:tree_m5}
\end{subfigure}
\hfill
\begin{subfigure}[c]{0.32\textwidth}
\centering
\begin{adjustbox}{max width=\linewidth}
\begin{tikzpicture}[scale=0.75, baseline=0]
    \tikzset{
        c1/.style={red, thick},
        c2/.style={blue, thick},
        c3/.style={green!70!black, thick},
        vnode/.style={circle, fill=black, inner sep=1.2pt}
    }

    \node[vnode, label={90:$\omega_0$}] (V0) at (90:2.1) {};
    \node[vnode, label={45:$\omega_1$}] (V1) at (45:2.1) {};
    \node[vnode, label={0:$\omega_2$}] (V2) at (0:2.1) {};
    \node[vnode, label={315:$\omega_3$}] (V3) at (315:2.1) {};
    \node[vnode, label={270:$\omega_4$}] (V4) at (270:2.1) {};
    \node[vnode, label={225:$\omega_5$}] (V5) at (225:2.1) {};
    \node[vnode, label={180:$\omega_6$}] (V6) at (180:2.1) {};
    \node[vnode, label={135:$\omega_7$}] (V7) at (135:2.1) {};

    \draw[c3] (V0) -- (V1); \draw[c3] (V2) -- (V3);
    \draw[c3] (V4) -- (V5); \draw[c3] (V6) -- (V7);

    \draw[c2] (V0) -- (V2); \draw[c2] (V0) -- (V3); \draw[c2] (V1) -- (V2); \draw[c2] (V1) -- (V3);
    \draw[c2] (V4) -- (V6); \draw[c2] (V4) -- (V7); \draw[c2] (V5) -- (V6); \draw[c2] (V5) -- (V7);

    \foreach \u in {V0, V1, V2, V3} {
        \foreach \v in {V4, V5, V6, V7} {
            \draw[c1] (\u) -- (\v);
        }
    }
\end{tikzpicture}
\end{adjustbox}
\caption{3-colored $K_8$}
\label{fig:k8_colored_5}
\end{subfigure}
\hfill
\begin{subfigure}[c]{0.32\textwidth}
\centering
\begin{adjustbox}{max width=\linewidth}
\begin{tikzpicture}[scale=0.75, baseline=0]
    \tikzset{
        c1/.style={red, thick},
        c2/.style={blue, thick},
        c3/.style={green!70!black, thick},
        vnode/.style={circle, fill=black, inner sep=1.2pt}
    }

    \node[vnode, label={90:$\omega_0$}] (P0) at (90:2.1) {};
    \node[vnode, label={45:$\omega_4$}] (P1) at (45:2.1) {};
    \node[vnode, label={0:$\omega_2$}] (P2) at (0:2.1) {};
    \node[vnode, label={315:$\omega_6$}] (P3) at (315:2.1) {};
    \node[vnode, label={270:$\omega_1$}] (P4) at (270:2.1) {};
    \node[vnode, label={225:$\omega_5$}] (P5) at (225:2.1) {};
    \node[vnode, label={180:$\omega_3$}] (P6) at (180:2.1) {};
    \node[vnode, label={135:$\omega_7$}] (P7) at (135:2.1) {};

    \draw[c1] (P0)--(P1); \draw[c1] (P1)--(P2); \draw[c1] (P2)--(P3); \draw[c1] (P3)--(P4);
    \draw[c1] (P4)--(P5); \draw[c1] (P5)--(P6); \draw[c1] (P6)--(P7); \draw[c1] (P7)--(P0);

    \draw[c1] (P0)--(P3); \draw[c1] (P1)--(P4); \draw[c1] (P2)--(P5); \draw[c1] (P3)--(P6);
    \draw[c1] (P4)--(P7); \draw[c1] (P5)--(P0); \draw[c1] (P6)--(P1); \draw[c1] (P7)--(P2);

    \draw[c2] (P0)--(P2); \draw[c2] (P1)--(P3); \draw[c2] (P2)--(P4); \draw[c2] (P3)--(P5);
    \draw[c2] (P4)--(P6); \draw[c2] (P5)--(P7); \draw[c2] (P6)--(P0); \draw[c2] (P7)--(P1);

    \draw[c3] (P0)--(P4);
    \draw[c3] (P1)--(P5);
    \draw[c3] (P2)--(P6);
    \draw[c3] (P3)--(P7);
\end{tikzpicture}
\end{adjustbox}
\caption{3-circulant representation}
\label{fig:circulant_m5}
\end{subfigure}

\caption{The binary tree $T$ of depth $m=3$ with branching numbers $k_1=k_2=k_3=2$, the corresponding 3-graph, and its circular representation $(1,3;2;4)_8$, obtained by ordering the leaves as $(\omega_0,\omega_4,\omega_2,\omega_6,\omega_1,\omega_5,\omega_3,\omega_7)$.}
\label{f5}
\end{figure}
}

\newcommand{\GallaiFigureSix}{%
\begin{figure}[htbp]
\centering
\begin{subfigure}[c]{\textwidth}
\centering
\begin{adjustbox}{max width=\linewidth}
\begin{tikzpicture}[x=0.82cm,y=0.82cm]
    \tikzset{
        tnode/.style={circle, draw=black, fill=white, inner sep=0.5pt,
            font=\scriptsize\sffamily, minimum size=13pt},
        leaf/.style={circle, fill=black, inner sep=1.2pt}
    }
    \foreach \i in {0,...,15} {
        \node[leaf, label={[font=\scriptsize]below:$\omega_{\i}$}]
            (L\i) at (\i,0) {};
    }
    \foreach \i in {0,...,7} {
        \pgfmathsetmacro{\x}{2*\i+0.5}
        \pgfmathtruncatemacro{\ca}{2*\i}
        \pgfmathtruncatemacro{\cb}{2*\i+1}
        \node[tnode] (N3_\i) at (\x,1) {4};
        \draw (N3_\i)--(L\ca); \draw (N3_\i)--(L\cb);
    }
    \foreach \i in {0,...,3} {
        \pgfmathsetmacro{\x}{4*\i+1.5}
        \pgfmathtruncatemacro{\ca}{2*\i}
        \pgfmathtruncatemacro{\cb}{2*\i+1}
        \node[tnode] (N2_\i) at (\x,2) {3};
        \draw (N2_\i)--(N3_\ca); \draw (N2_\i)--(N3_\cb);
    }
    \foreach \i in {0,1} {
        \pgfmathsetmacro{\x}{8*\i+3.5}
        \pgfmathtruncatemacro{\ca}{2*\i}
        \pgfmathtruncatemacro{\cb}{2*\i+1}
        \node[tnode] (N1_\i) at (\x,3) {2};
        \draw (N1_\i)--(N2_\ca); \draw (N1_\i)--(N2_\cb);
    }
    \node[tnode] (v0) at (7.5,4) {1};
    \draw (v0)--(N1_0); \draw (v0)--(N1_1);
    \node[circle, draw=black, fill=white, inner sep=0pt, minimum size=3pt,
        label={[font=\scriptsize]left:$s$}] (vs) at (7.5,4.65) {};
    \draw[densely dashed] (vs)--(v0);
\end{tikzpicture}
\end{adjustbox}
\caption{Game tree $T$ ($m=4$)}
\label{fig:tree_m6}
\end{subfigure}
\par\medskip
\begin{subfigure}[c]{0.48\textwidth}
\centering
\begin{adjustbox}{max width=\linewidth}
\begin{tikzpicture}[scale=1, baseline=0]
    \tikzset{
        c1/.style={red, line width=0.45pt},
        c2/.style={blue, line width=0.55pt},
        c3/.style={green!70!black, line width=0.7pt},
        c4/.style={magenta, line width=0.85pt},
        vnode/.style={circle, fill=black, inner sep=1.2pt}
    }

    \foreach \i in {0,...,15} {
        \pgfmathsetmacro{\angle}{90 - \i * 22.5}
        \node[vnode, label={[font=\scriptsize]\angle:$\omega_{\i}$}] (V\i) at (\angle:2.55) {};
    }

    \foreach \i in {0,...,7} {
        \foreach \j in {8,...,15} {
            \draw[c1] (V\i) -- (V\j);
        }
    }

    \foreach \b in {0,8} {
        \foreach \i in {0,1,2,3} {
            \foreach \j in {4,5,6,7} {
                \pgfmathtruncatemacro{\u}{\b+\i}
                \pgfmathtruncatemacro{\v}{\b+\j}
                \draw[c2] (V\u) -- (V\v);
            }
        }
    }

    \foreach \b in {0,4,8,12} {
        \foreach \i in {0,1} {
            \foreach \j in {2,3} {
                \pgfmathtruncatemacro{\u}{\b+\i}
                \pgfmathtruncatemacro{\v}{\b+\j}
                \draw[c3] (V\u) -- (V\v);
            }
        }
    }

    \foreach \i in {0,2,4,6,8,10,12,14} {
        \pgfmathtruncatemacro{\j}{\i+1}
        \draw[c4] (V\i) -- (V\j);
    }
\end{tikzpicture}
\end{adjustbox}
\caption{4-colored $K_{16}$}
\label{fig:k16_colored}
\end{subfigure}
\hfill
\begin{subfigure}[c]{0.48\textwidth}
\centering
\begin{adjustbox}{max width=\linewidth}
\begin{tikzpicture}[scale=1, baseline=0]
    \tikzset{
        c1/.style={red, line width=0.45pt},
        c2/.style={blue, line width=0.55pt},
        c3/.style={green!70!black, line width=0.7pt},
        c4/.style={magenta, line width=0.85pt},
        vnode/.style={circle, fill=black, inner sep=1.2pt}
    }

    \def\order{0, 8, 4, 12, 2, 10, 6, 14, 1, 9, 5, 13, 3, 11, 7, 15}
    \foreach \val [count=\i from 0] in \order {
        \pgfmathsetmacro{\angle}{90 - \i * 22.5}
        \node[vnode, label={[font=\scriptsize]\angle:$\omega_{\val}$}] (P\i) at (\angle:2.55) {};
    }

    \foreach \i in {0,...,14} {
        \pgfmathtruncatemacro{\start}{\i+1}
        \foreach \j in {\start,...,15} {
            \pgfmathtruncatemacro{\diff}{\j - \i}
            \pgfmathtruncatemacro{\dist}{\diff > 8 ? 16 - \diff : \diff}
            
            \ifnum\dist=1 \draw[c1] (P\i) -- (P\j); \fi
            \ifnum\dist=3 \draw[c1] (P\i) -- (P\j); \fi
            \ifnum\dist=5 \draw[c1] (P\i) -- (P\j); \fi
            \ifnum\dist=7 \draw[c1] (P\i) -- (P\j); \fi
            
            \ifnum\dist=2 \draw[c2] (P\i) -- (P\j); \fi
            \ifnum\dist=6 \draw[c2] (P\i) -- (P\j); \fi
            
            \ifnum\dist=4 \draw[c3] (P\i) -- (P\j); \fi
            
            \ifnum\dist=8 \draw[c4] (P\i) -- (P\j); \fi
        }
    }
\end{tikzpicture}
\end{adjustbox}
\caption{4-circulant representation}
\label{fig:circulant_m6}
\end{subfigure}

\caption{The binary tree $T$ of depth $m=4$ with branching numbers $k_1=k_2=k_3=k_4=2$, the corresponding 4-graph, and its circular representation $(1,3,5,7;2,6;4;8)_{16}$. In the circular order, the leaf indices are $(0,8,4,12,2,10,6,14,1,9,5,13,3,11,7,15)$. Colors 1, 2, 3, and 4 are red, blue, green, and magenta, respectively.}
\label{f6}
\end{figure}
}

\title{On Gallai $d$-colored circulants}
\author{Vladimir Gurvich\\
\texttt{vladimir.gurvich@gmail.com}\\
National Research University Higher School of Economics,\\
Moscow, Russian Federation;\\
RUTCOR, Rutgers University, Piscataway, NJ, United States
\and
Mariya Naumova\correspondingauthor\\
\texttt{mnaumova@business.rutgers.edu}\\
Rutgers Business School, Rutgers University,\\
Piscataway, NJ, United States}
\date{}

\begin{document}
\maketitle

\begin{abstract}
A $d$-graph is a complete graph whose edges are colored with colors from
$[d]$. It is Gallai if it contains no triangle whose three edges have
distinct colors. We prove that every Gallai $d$-circulant is an iterated
lexicographic product of circulants using at most two colors. The proof
uses the compatibility of the known canonical modular decomposition
with the cyclic action. For $d$-graphs avoiding both a rainbow triangle
$\Delta$ and the two-colored four-vertex graph $\Pi$ whose chromatic
components are paths, we give a more explicit classification.
Color-preserving vertex transitivity is equivalent to circular symmetry
and to a reduced positional game tree whose branching numbers and
labels depend only on depth. The edge colors in the given cyclic
labeling are determined by a chain of divisors of the order. This
description yields a complete isomorphism invariant and a linear-time
recognition procedure in the difference-word representation. The
characterization of CIS $d$-circulants reduces to CIS two-colored
factors. Finally, we extend the numerical CIS criterion to
$d$-graphs whose chromatic components are vertex-transitive.


\end{abstract}

\section{Introduction}

In this work, we assume all mentioned graphs are finite, simple, undirected, and have at least one vertex.
A $d$-\emph{graph} $G=(V;E_1,\ldots,E_d)$ is a complete graph
whose edge set is partitioned into $d$ color classes, some of which may
be empty. The graphs $G_i=(V,E_i)$, $i\in[d]=\{1,\ldots,d\}$, are called
the \emph{chromatic components} of $G$. We distinguish the number $d$
of available colors from the number of nonempty color classes. A
two-colored graph is an ordinary graph together with its complement.

We write $\varphi_G(xy)$ for the color of the edge joining distinct
vertices $x,y$. All isomorphisms and automorphisms preserve each color
individually. When discussing a forbidden colored subgraph, however,
its distinct colors may be renamed. All forbidden subgraphs are induced
on their vertex sets.

Two small colored graphs will play an important role. The graph
$\Delta$ is a triangle with three pairwise distinct edge colors. The
graph $\Pi$ has vertices $v_1,v_2,v_3,v_4$ and two color classes
\[
E_1=\{v_1v_2,v_2v_3,v_3v_4\},\qquad
E_2=\{v_2v_4,v_4v_1,v_1v_3\}.
\]
Thus both chromatic components of $\Pi$ are isomorphic to $P_4$, the
path on four vertices. We regard these patterns as $d$-graphs by adding
empty color classes when necessary; see Fig.~\ref{f1}.

\GallaiFigureOne

A $\Delta$-free $d$-graph is called a \emph{Gallai $d$-graph}. Gallai's
decomposition theorem and its later formulations describe such graphs
by substitutions of graphs using at most two colors
\cite{Gal67,GS04,Gur09}; see also \cite{CE97,CEL86,CG83,GN26}.
The $d$-graphs avoiding both $\Pi$ and $\Delta$ admit a more restricted
canonical decomposition. They correspond to reduced positional game
structures on trees \cite{Gur78,Gur82,Gur84,Gur09}. For two colors this
is the familiar cotree representation of $P_4$-free graphs, closely
related to monotone read-once Boolean functions
\cite{CLB81,GG11,Gur77,Gur91,Sei74,Sum73}.

The same labeled-tree representation is studied in the literature on
\emph{symbolic ultrametrics} \cite{BD98,HW18}. After a special symbol is
assigned to the diagonal, the conditions defining a symbolic
ultrametric are precisely the exclusion of $\Delta$ and $\Pi$.
We shall use the game-theoretic terminology, while keeping this
equivalent interpretation in mind.

In this paper we determine the restrictions imposed by circular
symmetry on these canonical decompositions.
Theorem~\ref{thm:gallai-circulants} shows that a
Gallai $d$-circulant is an iterated lexicographic product of circulants
using at most two colors. The essential observation is that the cyclic
action preserves the canonical strong-module partition. It follows
that the quotient and every block are circulants and that all blocks
are color-isomorphic. For CIS $d$-circulants, the resulting
characterization reduces the problem to CIS two-colored factors.

For the $(\Pi,\Delta)$-free class, Theorem~\ref{thm:uniform-tree}
shows that color-preserving vertex transitivity is equivalent to
circular symmetry and to a reduced tree whose branching numbers and
player labels depend only on depth. All canonical quotient factors
are then monochromatic. Theorem~\ref{thm:divisibility} describes the
edge colors by a chain of divisors in the given cyclic labeling.
The resulting parameter sequence is a complete isomorphism invariant
and can be recovered in linear time from the difference word; see
Corollary~\ref{cor:recognition}. We also extend the numerical CIS
criterion to $d$-graphs whose chromatic components are vertex-transitive.

The general Gallai decomposition, the general CIS decomposition into
factors using at most two colors, and the two-color case of the
$P_4$-free circulant decomposition are established results
\cite{GS04,LT26,BGM14}. The additional conclusions here are the
preservation of circular symmetry by the canonical factors, the
vertex-transitive classification in the $(\Pi,\Delta)$-free case,
and its explicit arithmetic and game-tree descriptions.

\section{CIS graphs and \texorpdfstring{$d$}{d}-graphs}
\label{sec:cis}

Given a $d$-graph $G=(V;E_1,\ldots,E_d)$, choose a maximal stable set
$S_i$ in each chromatic component $G_i$, and put
\[
S=\bigcap_{i=1}^d S_i.
\]
Here and below, maximality is with respect to inclusion. We have
$|S|\le1$: if two distinct vertices belonged to every $S_i$, the edge
between them could have none of the $d$ colors.

We say that $G$ has the \emph{CIS property}, or is a \emph{CIS
$d$-graph}, if $S\ne\emptyset$ for every such choice. In that case
$|S|=1$. An empty color class causes no difficulty, since the only
maximal stable set of its chromatic component is $V$.
For $d=2$, the condition says that every maximal clique meets every
maximal stable set of either of the two complementary graphs. This
explains the name CIS. For background on the ordinary and multicolored
properties, see \cite{ABG06,ABG10,BG09,BGM17,Gur11}.
Further characterizations and questions for ordinary CIS graphs are
studied in \cite{DLZ04,DLZ04a,Zan95}.

The assertion that every CIS $d$-graph is Gallai was known as the
$\Delta$-conjecture. It was formulated in
\cite[remark following Claim 17, p.~71]{Gur78}. Earlier work established
that every $(\Pi,\Delta)$-free $d$-graph is CIS and that it suffices to
prove the conjecture for three colors. The latter reduction was
attributed to Andrey Gol'berg in 1975; see \cite{Gur09,GN26,LT26} for
the history and related results.

\begin{theorem}[Liu and Tang \cite{LT26}]
\label{thm:delta}
Every CIS $d$-graph is Gallai.
\end{theorem}

The same paper records the resulting canonical decomposition of CIS
$d$-graphs into CIS factors using at most two colors
\cite[Corollary 4.3]{LT26}. We recall this consequence in
Section~\ref{sec:gallai} before applying it to circulants.

Also in 2026, Tao, Yang, and Zang proved that recognizing ordinary CIS
graphs is coNP-complete \cite{TYZ26}. This result concerns arbitrary
graphs; it does not by itself establish the complexity of recognition
within the class of circulants.

\begin{remark}
An ordinary graph is called \emph{almost CIS} if exactly one pair
consisting of a maximal clique and a maximal stable set is disjoint.
The relation between almost CIS graphs and split graphs was studied in
\cite{BGZ09}. Wu, Zang, and Zhang proved that a graph is almost CIS if
and only if it has a unique split partition \cite{WZZ09}.
\end{remark}

\section{Substitution and lexicographic product}
\label{s-sub-lp}

For vertex-disjoint $d$-graphs $G,F$ and $v\in V(G)$, substituting $F$
for $v$ means deleting $v$ and joining each $x\in V(F)$ to each
$u\in V(G)\setminus\{v\}$ by an edge of color $\varphi_G(vu)$.
The colors within $F$ and within $G-v$ remain unchanged. Denote the
result by $G_v[F]$. If the palettes differ, we first embed them in their
union by adding empty color classes.

For two colors this is ordinary graph substitution, with edges and
nonedges receiving the two colors. Families preserved by substitution
occur throughout graph theory. A classical example is the family of
perfect graphs; the substitution lemma is closely connected with the
weak perfect graph theorem \cite{Ful71,Lov72a,Lov72b}.

We call a family $\Phi$ \emph{exactly closed under substitution} if
\[
G_v[F]\in\Phi\quad\Longleftrightarrow\quad G,F\in\Phi.
\]
Only nonempty graphs are substituted.

\begin{lemma}
\label{lem:substitution}
The families of $\Delta$-free $d$-graphs, $\Pi$-free $d$-graphs, and
CIS $d$-graphs are exactly closed under substitution.
\end{lemma}

\begin{proof}
Neither $\Pi$ nor $\Delta$ has a proper module with at least two
vertices. Here a module is a set whose vertices have the same colored
adjacencies to every vertex outside the set. In $G_v[F]$, the vertices
of $F$ form a module. A forbidden copy meeting both $V(F)$ and its
complement can therefore use only one vertex of $F$; replacing
that vertex by $v$ gives the same forbidden copy in $G$. A copy entirely
inside $F$ is already a copy in $F$. Conversely, both $G$ and $F$ occur
as induced colored subgraphs of $G_v[F]$. This proves the first two
assertions.

For the CIS assertion, let $S_i$ be a maximal stable set in the $i$th
chromatic component of $G_v[F]$. If $S_i$ meets $F$, replace its
intersection with $F$ by $v$; otherwise leave it unchanged. The
resulting set $U_i$ is maximal stable in $G_i$. In the first case
$S_i\cap V(F)$ is maximal stable in $F_i$. These facts follow directly
from maximality and the definition of substitution.

Suppose that both $G$ and $F$ are CIS. The sets $U_i$ have a common
vertex. If it is different from $v$, it belongs to all $S_i$. If it is
$v$, all $S_i$ meet $F$, and their restrictions to $F$ have a common
vertex by the CIS property of $F$.

Conversely, maximal stable sets in $G_i$ with empty common intersection
can be lifted to $G_v[F]$ by replacing $v$, whenever present, by a
maximal stable set in $F_i$. The lifted sets are maximal and still have
empty common intersection. If instead $F$ is not CIS, choose maximal
stable sets in its components with empty intersection, and choose in
each $G_i$ a maximal stable set containing $v$. Their common
intersection in $G$ is exactly $\{v\}$. Substituting the chosen sets
from $F$ again gives maximal stable sets with empty intersection.
\end{proof}

The CIS part of Lemma~\ref{lem:substitution} is a standard property of
these graphs; see \cite{ABG06,Gur09,BGM14}. Notice that it does not
assert that the CIS property is inherited by arbitrary induced
subgraphs.

A $d$-graph is \emph{complementary connected} (CC) if the complement
of each chromatic component is connected. When $|V(G)|\ge2$,
\[
G_v[F]\text{ is CC}\quad\Longleftrightarrow\quad G\text{ is CC}.
\]
Indeed, contracting the substituted set preserves the necessary
connectivity, while a connected complement of $G$ has a neighbor of
$v$ outside that set, joining all its substituted vertices. Thus CC
graphs are closed under substitution, but are not exactly closed.
The restriction $|V(G)|\ge2$ excludes the exceptional case in which
the substitution simply gives $F$.

Substituting the same graph $F$ for every vertex of $G$ gives the
\emph{lexicographic product}, denoted by $G[F]$ or $G\times F$.
Its vertex set is $V(G)\times V(F)$, and
\begin{equation}
\varphi_{G[F]}\bigl((u,x)(v,y)\bigr)=
\begin{cases}
\varphi_G(uv),&u\ne v,\\
\varphi_F(xy),&u=v,
\end{cases}
\label{eq:lex-definition}
\end{equation}
for distinct vertices. This operation is associative up to the natural
color-preserving isomorphism, but is not commutative. If both factors
have at least two vertices, its used colors are exactly the union of
their used colors. We use the nested notation $G^1[G^2[\cdots[G^s]\cdots]]$
when the order of the factors is important.

\section{Circulants and \texorpdfstring{$d$}{d}-circulants}
\label{sec:circulants}

A \emph{$d$-circulant} of order $n$ is a $d$-graph on $\Z_n$ whose
coloring has the form
\[
\varphi_G(xy)=c(y-x),\qquad
c:\Z_n\setminus\{0\}\longrightarrow[d],\qquad c(a)=c(-a).
\]
Equivalently, the shift $x\mapsto x+1$ is a color-preserving
automorphism. We also call a $d$-graph a circulant if it admits such a
labeling. One may specify $c$ by partitioning the circular edge
lengths $1,\ldots,\lfloor n/2\rfloor$ into color classes. For example,
\[
(1,3,5,7;2,6;4;8)_{16}
\]
has four colors, listed in order and separated by semicolons. For two
colors we obtain a pair of complementary ordinary circulants.

There is a direct arithmetic test for the Gallai property:
\begin{equation}
G\text{ is Gallai}\quad\Longleftrightarrow\quad
\bigl|\{c(a),c(b),c(a-b)\}\bigr|\le2
\label{eq:gallai-arithmetic}
\end{equation}
for all distinct nonzero $a,b\in\Z_n$. To prove it, translate one
vertex of a triangle to $0$. The remaining vertices are $a,b$, and the
three edge colors are those in \eqref{eq:gallai-arithmetic}.

Every two-colored circulant is Gallai. A three-colored circulant of
even order is also Gallai if one color class consists of all odd edge
lengths. Every triangle then contains either two odd-length edges,
which have the same color, or no odd-length edge, in which case at
most the other two colors occur. In particular, $(1,3;2;4)_8$ is
Gallai.

For an ordinary circulant $C_n(L)$ with generating lengths $L$, the
number of connected components is
\[
g=\gcd\bigl(n,\{\ell:\ell\in L\}\bigr).
\]
Each component is isomorphic to $C_{n/g}(L/g)$. For $L\ne\emptyset$,
the exact bipartiteness criterion is
\begin{equation}
C_n(L)\text{ is bipartite}\quad\Longleftrightarrow\quad
n/g\text{ is even and every }\ell/g\text{ is odd}.
\label{eq:bipartite}
\end{equation}
For a connected circulant this reduces to the familiar criterion that
$n$ is even and all generators are odd
\cite{BGM14}. An edgeless graph is bipartite. The
connectedness qualification cannot be omitted: $C_8(\{2\})$ is the
disjoint union of two four-cycles.

Recall that an ordinary graph is \emph{well-covered} if all its
maximal stable sets have the same size, necessarily its independence
number $\alpha$. For every ordinary circulant $H$,
\begin{equation}
\alpha(H)\omega(H)\le |V(H)|.
\label{eq:alpha-omega}
\end{equation}
Moreover, $H$ is CIS if and only if equality holds and both $H$ and
its complement are well-covered \cite{BGM14}. This
criterion also holds for vertex-transitive graphs \cite{DHMV15}.
Its multicolored form is given in Section~\ref{sec:numerical}.

\begin{proposition}
\label{prop:circulant-product}
The lexicographic product of two colored circulants is a colored
circulant. If $G$ has order $n$ and $F$ has order $m$, a difference
coloring of $G[F]$ on $\Z_{nm}$ is
\begin{equation}
c_{G[F]}(a)=
\begin{cases}
c_G(a\bmod n),& n\nmid a,\\[1mm]
c_F(a/n),&n\mid a,
\end{cases}
\qquad 1\le a<nm.
\label{eq:product-colors}
\end{equation}
\end{proposition}

\begin{proof}
Identify $(u,x)\in\Z_n\times\Z_m$ with $u+nx\in\Z_{nm}$, using
$0\le u<n$ and $0\le x<m$. Two vertices are in the same substituted
copy of $F$ exactly when their difference is divisible by $n$.
Otherwise their edge has the color of the corresponding difference
modulo $n$ in $G$. This gives \eqref{eq:product-colors}, which is
symmetric under $a\mapsto nm-a$.
\end{proof}

For ordinary circulants this is the usual product formula
\cite{BGM14}; see also \cite{GGSTU93,Hos08}.
The order of the factors matters, and their colors are retained.
If $G$ is a red $K_2$ and $F$ a blue
$K_3$, then $G[F]=(1,3;2)_6$, whereas $F[G]=(3;1,2)_6$.

\begin{remark}
Circulants provide useful examples in the study of graph classes and
conjectures. They have been used for perfect and partitionable graphs
\cite{BG93,GGSTU93,BBGMP98}, graphs without even holes \cite{ABG05},
and oriented graphs without kernels \cite{AFG98,GN26}.
The arithmetic description of their edges often makes constructions
and computations particularly convenient.
\end{remark}

\section{Decomposing Gallai and CIS \texorpdfstring{$d$}{d}-circulants}
\label{sec:gallai}

We first recall the part of modular decomposition needed below. A set
$M\subseteq V(G)$ is a \emph{module} if, for each $x\notin M$, the
color $\varphi_G(xy)$ is independent of $y\in M$. Two sets overlap if
their intersection and both of their differences are nonempty. A
module is \emph{strong} if it overlaps no module. For a graph with at
least two vertices, the maximal proper strong modules form a canonical
partition of its vertex set. The quotient has one vertex for each
part, with the color between parts given by the color of all edges
between them.

Recursively applying this construction gives the canonical modular
decomposition tree. Its leaves are the original vertices, and its
internal nodes carry quotient graphs. For symmetric colored complete
graphs each such quotient is monochromatic or prime, where prime
means that the only nonempty modules are singleton sets and the whole
vertex set.
We use these standard facts about modular decomposition
\cite{Mor85,McS99,MS89,ER90,Gur09}; their colored formulation is recalled in
\cite[Corollary 4.3 and its proof]{LT26}.

Each quotient is color-isomorphic to the induced subgraph obtained
by choosing one representative from each part. In particular, a
quotient of a Gallai $d$-graph is again Gallai.

Gallai's theorem states that every Gallai $d$-graph of order at least
two has a partition into at least two nonempty modules such that at
most two colors occur between the parts
\cite{Gal67,GS04}. In a prime quotient all parts of such a partition
must be singletons. Therefore every canonical quotient of a Gallai
$d$-graph uses at most two colors. Conversely, substitutions of
one- and two-colored graphs remain Gallai by
Lemma~\ref{lem:substitution}.

Together with Theorem~\ref{thm:delta} and exact CIS closure, this
implies that every CIS $d$-graph has a canonical modular decomposition
into CIS factors using at most two colors
\cite[Corollary 4.3]{LT26}. The following observation gives the extra
structure supplied by circular symmetry.

\begin{lemma}
\label{lem:cyclic-modules}
Let $G$ be a colored circulant on $\Z_n$, $n\ge2$, and let
$B_0,\ldots,B_{q-1}$ be its maximal proper strong modules, with
$0\in B_0$. Then $n=qm$ for some $m$, and the blocks can be indexed as
\[
B_r=\{r,r+q,\ldots,r+(m-1)q\},\qquad 0\le r<q.
\]
The quotient $Q$ is a circulant of order $q$, the induced graph
$F=G[B_0]$ is a circulant of order $m$, and $G\cong Q[F]$.
\end{lemma}

\begin{proof}
The shift $\tau(x)=x+1$ preserves the canonical partition. Its cyclic
group is transitive on the blocks, so they have the same size. The
stabilizer of $B_0$ is a subgroup of $\Z_n$. It is exactly $B_0$:
if $r\in B_0$, then $\tau^r(B_0)$ meets $B_0$ at $r$ and hence is the
same block; the converse follows by applying $\tau^r$ to $0$.
Thus $B_0=q\Z_n$ and the other blocks are its cosets.

The shift induces a $q$-cycle on the quotient, and $\tau^q$ induces an
$m$-cycle on each block. Hence $Q$ and $F$ are circulants. The shifts
$\tau^r$ identify all the induced blocks with $F$. Since the blocks
are modules, their mutual adjacencies are exactly those prescribed by
$Q$, giving $G\cong Q[F]$.
\end{proof}

\begin{theorem}[Gallai circulant factorization]
\label{thm:gallai-circulants}
A colored circulant $G$ is Gallai if and only if it admits a
representation
\begin{equation}
G\cong Q_1[Q_2[\cdots[Q_s]\cdots]],
\label{eq:gallai-factorization}
\end{equation}
where each $Q_j$ is a circulant using at most two colors. For any such
representation, $G$ is CIS if and only if every $Q_j$ is CIS.
\end{theorem}

\begin{proof}
For $|V(G)|=1$ use the empty product, interpreted as the one-vertex
graph. Otherwise, apply Lemma~\ref{lem:cyclic-modules}. The quotient
$Q$ uses at most two colors by the Gallai decomposition just recalled.
The induced block $F$ is again a Gallai circulant, and has fewer
vertices. Induction gives \eqref{eq:gallai-factorization}, omitting
one-vertex factors.

Conversely, every one- or two-colored graph is Gallai. Substitution
preserves this property by Lemma~\ref{lem:substitution}, and the
product is a circulant by Proposition~\ref{prop:circulant-product}.
The CIS assertion follows by repeated application of the same
substitution lemma.
\end{proof}

Thus the canonical tree does not become a path: it has isomorphic
subtrees at every fixed depth. Recording one quotient at each depth
gives the sequence of factors in \eqref{eq:gallai-factorization}.
The colorings of the factors determine the product coloring by
\eqref{eq:product-colors}.
Although the representations in Theorem~\ref{thm:gallai-circulants}
need not be unique, the canonical modular decomposition determines
a unique ordered sequence of quotient isomorphism types. For example,
$K_4^{\,1}\cong K_2^{\,1}[K_2^{\,1}]$, but only the single factor
$K_4^{\,1}$ is canonical.

One-colored factors are complete graphs and are automatically CIS.
The characterization of CIS $d$-circulants therefore reduces to CIS
two-colored circulants. Large families of the latter are studied in
\cite{BGM14}. This reduction should be distinguished from a complete
constructive description of all the two-colored building blocks.

\begin{corollary}
\label{cor:prime-gallai}
A Gallai circulant of prime order uses at most two colors. More
generally, if all factors in \eqref{eq:gallai-factorization} have
order at least two, their number is at most $\Omega(n)$, where
$n=|V(G)|$ and $\Omega(n)$ counts prime factors with multiplicity.
If $G$ is CIS and uses $r$ colors, then $r\le\Omega(n)$.
In particular, a CIS circulant of prime order is monochromatic.
\end{corollary}

\begin{proof}
The factor orders multiply to $n$. Each nontrivial factor contributes
at least one prime factor.
For the CIS assertion, omit one-vertex factors and let $q_j$ and
$r_j\in\{1,2\}$ be the order and the number of used colors of $Q_j$.
If $r_j=1$, then $\Omega(q_j)\ge1$. If $r_j=2$, let $H_j$ be one of
the two nonempty chromatic components of this factor. Then $H_j$ is
an ordinary CIS circulant with
$\alpha(H_j)\omega(H_j)=q_j$ and
$\alpha(H_j),\omega(H_j)\ge2$. Hence $\Omega(q_j)\ge2$.
The colors used by the product form the union of those used by its
factors, so
\[
r\le\sum_j r_j\le\sum_j\Omega(q_j)=\Omega(n).\qedhere
\]
\end{proof}

\section{Decomposing \texorpdfstring{$\Pi$- and $\Delta$-free $d$}{Pi- and Delta-free d}-circulants}
\label{sec:trees}

Consider a finite rooted tree $T$ with root $v_0$, in which each
internal vertex has at least two children. Assign a label
$i(v)\in[d]$ to every internal vertex, with distinct labels at
adjacent internal vertices. We call such a labeled tree \emph{reduced}.
It defines a positional game structure: player $i(v)$ chooses an
outgoing edge at position $v$. The game starts at $v_0$ and terminates
upon reaching a leaf. There are no forced moves, and no player moves
twice successively along a play.

In the figures we retain an auxiliary uncolored vertex $s$ joined to
$v_0$ by a dashed stem. This vertex only indicates the root. It is not
a game position or an outcome, and the stem is not a move. The leaves
of $T$ below $v_0$, excluding $s$, are the terminal outcomes.

For two distinct leaves $t,t'$, their last common position is their
lowest common ancestor $\lca_T(t,t')$. Define the complete colored
graph $G(T)$ on the leaves by
\begin{equation}
\varphi_{G(T)}(tt')=i\bigl(\lca_T(t,t')\bigr).
\label{eq:lca-color}
\end{equation}

The following known correspondence is the starting point for this
section \cite[Chapter 5]{Gur78}; see also \cite{Gur82,Gur84,Gur09}.
Its symbolic-ultrametric formulation is given in \cite{BD98,HW18}.

\begin{theorem}[Reduced tree representation]
\label{thm:tree-known}
A complete colored graph is $(\Pi,\Delta)$-free if and only if it is
$G(T)$ for a reduced rooted labeled tree $T$. The reduced tree is
unique up to a rooted label-preserving isomorphism that respects the
identification of its leaves with the graph vertices.
\end{theorem}

The singleton graph corresponds to the tree with one leaf and no
internal vertices. For larger graphs, the internal quotients are
monochromatic complete graphs. Conversely, substituting these graphs
gives \eqref{eq:lca-color}. The reduction conditions remove unary
vertices and contract adjacent internal vertices with the same label.
This is the canonical tree in Theorem~\ref{thm:tree-known}.

For two colors, interpreting the internal labels as conjunction and
disjunction yields the usual representation of monotone read-once
Boolean functions \cite{Gur77,Gur91,GG11}. Arbitrary read-once
functions may also have negated literals, which require additional
leaf signs. Related Boolean and hypergraph formulations are discussed
in \cite{Eit94,Eit95,Gur82,Gur84}.

We now transfer symmetry from the graph to its canonical tree.
A $d$-graph is \emph{vertex-transitive} if its color-preserving
automorphism group is transitive on vertices.

\begin{theorem}[Classification by uniform reduced trees]
\label{thm:uniform-tree}
Let $G$ be a $(\Pi,\Delta)$-free $d$-graph with $n\ge2$ vertices.
The following properties are equivalent: $G$ is vertex-transitive;
$G$ is a circulant; and its reduced tree has all leaves at one depth
$m$, with a common label $i_j$ and a common number $k_{j+1}$ of
children at each internal depth $j$, $0\le j<m$.
In this case
\begin{equation}
G\cong K_{k_1}^{\,i_0}
 \bigl[K_{k_2}^{\,i_1}[\cdots[K_{k_m}^{\,i_{m-1}}]\cdots]\bigr],
\label{eq:mono-product}
\end{equation}
where $K_q^{\,i}$ is the complete graph of order $q$ with all edges of
color $i$, $k_j\ge2$, and $i_j\ne i_{j+1}$.
\end{theorem}

\begin{proof}
Every color-preserving automorphism of $G$ permutes the leaves of its
reduced tree. By uniqueness in Theorem~\ref{thm:tree-known}, this
permutation extends to a root- and label-preserving tree automorphism.
If $G$ is vertex-transitive, all leaves consequently have the same
depth $m$.

Choose two vertices $u,v$ at the same depth and a leaf below each of
them. An automorphism taking the first leaf to the second takes its
ancestor at that depth to the corresponding ancestor of the second
leaf; hence it takes $u$ to $v$. Thus the labels and numbers of
children are constant at each depth. The reduction conditions give
$k_j\ge2$ and $i_j\ne i_{j+1}$, and the lowest-common-ancestor rule
gives \eqref{eq:mono-product}.

Conversely, every monochromatic complete graph is a circulant, so
\eqref{eq:mono-product} is a circulant by
Proposition~\ref{prop:circulant-product}. Every circulant is
vertex-transitive.
\end{proof}

The theorem uses vertex transitivity rather than a prescribed cyclic
action. In particular, every vertex-transitive $(\Pi,\Delta)$-free
$d$-graph is a circulant. For two colors, the circulant decomposition
in \eqref{eq:mono-product} specializes to the known decomposition of
$P_4$-free circulants \cite[Section 5]{BGM14}.

The tree description has a simple arithmetic form. Put
\begin{equation}
P_0=1,\qquad P_j=\prod_{h=1}^{j}k_h\quad(1\le j\le m),
\qquad P_m=n.
\label{eq:prefix-products}
\end{equation}

\begin{theorem}[Divisor-chain characterization]
\label{thm:divisibility}
A $d$-circulant on $\Z_n$, $n\ge2$, is $(\Pi,\Delta)$-free if and
only if there are integers $k_1,\ldots,k_m\ge2$ with product $n$ and
labels $i_0,\ldots,i_{m-1}\in[d]$, distinct at consecutive indices,
such that its difference coloring satisfies
\begin{equation}
c(a)=i_j\quad\text{when}\quad
P_j\mid a\ \text{ and }\ P_{j+1}\nmid a,
\qquad 1\le a<n.
\label{eq:divisibility-colors}
\end{equation}
The numbers $P_j$ are defined in \eqref{eq:prefix-products}; the rule
holds for the given labeling on $\Z_n$.
\end{theorem}

\begin{proof}
Let $T$ be the uniform reduced tree of the circulant. The shift
$\tau(x)=x+1$ lifts to an automorphism of $T$. For each depth $j$,
the leaf sets below the vertices at that depth form a
$\tau$-invariant partition with $P_j$ parts. As in the proof of
Lemma~\ref{lem:cyclic-modules}, these parts are exactly the residue
classes modulo $P_j$.

The leaves $0,a$ have the same ancestor at depth $j$ precisely when
$P_j\mid a$. Their lowest common ancestor is at depth $j$ when they
first separate at depth $j+1$. This proves
\eqref{eq:divisibility-colors} for the given cyclic labeling.

Conversely, make the uniform tree with the specified branching numbers
and labels. Label a leaf by its successive child indices
$(a_1,\ldots,a_m)$, where $0\le a_h<k_h$, and identify it with
\[
a_1+k_1a_2+k_1k_2a_3+\cdots+P_{m-1}a_m\in\Z_n.
\]
Two leaves agree through depth $j$ exactly when their difference is
divisible by $P_j$. Thus their lowest-common-ancestor coloring is
\eqref{eq:divisibility-colors}. It is symmetric because every $P_j$
divides $n$. Theorem~\ref{thm:tree-known} now gives the assertion.
\end{proof}

\begin{corollary}
\label{cor:recognition}
The ordered sequence $(k_1,i_0),\ldots,(k_m,i_{m-1})$ is a complete
invariant for color-preserving isomorphism among vertex-transitive
$(\Pi,\Delta)$-free $d$-graphs of order at least two. Given the
explicit difference word $c(1),\ldots,c(n-1)$ of a $d$-circulant,
one can recognize $(\Pi,\Delta)$-freeness and recover this sequence
in $O(n)$ time.
\end{corollary}

\begin{proof}
The invariant assertion follows from uniqueness of the reduced tree
and Theorem~\ref{thm:uniform-tree}. For recognition, the singleton
graph gives the empty sequence. Otherwise, put $i=c(1)$.
If all entries equal $i$, the graph is $K_n^{\,i}$. Otherwise, put
\[
q=\min\{1\le a<n:c(a)\ne i\}.
\]
In a reduced representation, $q=k_1$: all smaller positive differences
have color $i_0$, whereas $c(k_1)=i_1\ne i_0$. Check that $q\mid n$
and that $c(a)=i$ whenever $q\nmid a$, and then apply the same
procedure to $c'(b)=c(qb)$, $1\le b<n/q$. Reject if any check fails.
Successful checks give $K_q^{\,i}[F]$, where $F$ has difference word
$c'$, and recursive application gives exactly
\eqref{eq:divisibility-colors}. Each recursive order is at most half
its predecessor, so the total work is
$O(n+n/2+n/4+\cdots)=O(n)$ in the explicit difference-word
representation, with unit-cost color comparisons and integer operations.
\end{proof}

The level labels need not be pairwise distinct; only successive
labels differ. If $r$ colors are
used, then
\begin{equation}
r\le m\le\Omega(n).
\label{eq:color-bound}
\end{equation}
In particular, a $(\Pi,\Delta)$-free circulant of prime order is
monochromatic. Such a graph is always CIS, since every factor in
\eqref{eq:mono-product} is CIS and
Lemma~\ref{lem:substitution} applies.

For two distinct leaves, their path within $T$ has even length. If the length
is $2h$, their lowest common ancestor is at depth $m-h$, so their
edge has color $i_{m-h}$. This gives the path-length interpretation of
the same coloring.

Finally, the players, or colors, may be partitioned into groups, with
each group replaced by one player. Apply this change to all internal
labels of $T$. If adjacent labels become equal, contract the edge
between them. Each contraction preserves the leaf coloring and leaves
at least two children at every internal vertex. The resulting reduced
tree shows that color merging preserves $(\Pi,\Delta)$-freeness and
hence the CIS property. Circular symmetry is preserved directly by
recoloring. This argument explains the preservation within the present
class; color merging can create a $\Pi$ in an arbitrary colored graph.

\section{A numerical CIS criterion}
\label{sec:numerical}

We extend the numerical criterion for ordinary vertex-transitive
CIS graphs to $d$-graphs. In fact, it is enough that each chromatic
component be vertex-transitive. For two colors it reduces to
\cite{BGM14} for circulants and
\cite[Theorem 3.1]{DHMV15} for vertex-transitive graphs.

\begin{proposition}[Multicolored numerical CIS criterion]
\label{prop:numerical-cis}
Let $G=(V;E_1,\ldots,E_d)$ be a $d$-graph of order $n\ge1$ whose
chromatic components $G_i$ are vertex-transitive. Then
\begin{equation}
\prod_{i=1}^d\alpha(G_i)\le n^{d-1}.
\label{eq:multi-inequality}
\end{equation}
Moreover, $G$ is CIS if and only if every $G_i$ is well-covered and
equality holds in \eqref{eq:multi-inequality}.
\end{proposition}

\begin{proof}
For each $i$, let $A_i=\Aut(G_i)$. Choose stable sets $S_i$ in $G_i$
and independent uniformly distributed $\gamma_i\in A_i$.
Each $\gamma_i(S_i)$ is stable in $G_i$, so their intersection has at
most one vertex. Transitivity of each $A_i$ gives
\begin{equation}
\E\left|\bigcap_{i=1}^d\gamma_i(S_i)\right|
=\sum_{v\in V}\prod_{i=1}^d\frac{|S_i|}{n}
=\frac{\prod_{i=1}^d|S_i|}{n^{d-1}}.
\label{eq:intersection-average}
\end{equation}
Choosing maximum stable sets proves \eqref{eq:multi-inequality}.

Suppose that $G$ is CIS and that the $S_i$ are maximal. Every
intersection on the left of \eqref{eq:intersection-average} then has
size one. Consequently $\prod_i|S_i|=n^{d-1}$ for every such choice.
Fixing all sets except $S_i$ shows that all maximal stable sets of
$G_i$ have the same size. Thus every component is well-covered, and
the required equality follows.

Conversely, suppose that all components are well-covered and that
equality holds. For any choice of maximal stable sets, the expectation
in \eqref{eq:intersection-average} equals one. The intersection size
takes only the values zero and one, so it equals one for every tuple
of automorphisms. In particular it equals one when all $\gamma_i$
are the identity. Hence $G$ is CIS.
\end{proof}

An unused color has $\alpha(G_i)=n$ and multiplies both sides of the
equality by $n$, so the statement is independent of unused colors.
In particular, the proposition applies when the color-preserving
automorphism group of $G$ is transitive. It is a numerical
characterization; its use requires the independence numbers and
well-coveredness of the chromatic components.
The well-coveredness condition is essential. For example, the
ordinary six-cycle satisfies $\alpha\omega=6$ but has a maximal
stable set $\{2,5\}$ disjoint from the maximal clique $\{0,1\}$.

For the uniform trees of Section~\ref{sec:trees}, the criterion can be
read directly from the branching numbers. Within any subtree that it
meets, a maximal stable set in the component of color $i$ chooses one
child at a vertex labeled $i$ and meets every child at a vertex
carrying another label. It follows by
induction that this component is well-covered and
\[
\alpha(G_i)=\prod_{\substack{0\le j<m\\ i_j\ne i}} k_{j+1}.
\]
Multiplying over $i\in[d]$ gives $n^{d-1}$, in agreement with
Proposition~\ref{prop:numerical-cis}.
For example, in $(1,3;2;4)_8$ all three components are well-covered
with independence number $4$, and $4^3=8^2$.

\section{Examples and game structures}
\label{sec:examples}

Table~\ref{tab:examples} distinguishes the three properties for several
small circulants.

\begin{table}[htbp]
\centering
\small
\begin{tabular}{@{}lccc@{}}
\toprule
Circulant & Gallai & $\Pi$-free & CIS\\
\midrule
$(1,3;2;4)_8$ & Yes & Yes & Yes\\
$(1,3,5;2;4)_{10}$ & Yes & No & No\\
$(1,3,5;2,6;4)_{12}$ & Yes & Yes & Yes\\
$(1,3,5,7;2,6;4)_{14}$ & Yes & No & No\\
$(1,3,5,7,9;2,6,8;4)_{18}$ & Yes & No & No\\
$(1,3,5,7;2,6;4;8)_{16}$ & Yes & Yes & Yes\\
$H_{36}$ & Yes & No & Yes\\
\bottomrule
\end{tabular}
\caption{Gallai, $\Pi$-free, and CIS properties of the listed circulants;
$H_{36}$ is defined below.}
\label{tab:examples}
\end{table}

All the three-colored examples have even order and a color class
consisting of all odd lengths, so they are Gallai. For orders
$10,14,18$, write the coloring as $K_2^{\,1}[F]$ using the even and odd
residue classes. In the respective inner graphs of orders $5,7,9$,
color $3$ gives a cycle (generated by difference $2$), while color $2$
gives its complement. Each of these cycles contains an induced
$P_4$ and is not CIS: for an odd cycle $C_q$, $q\ge5$,
$\alpha(C_q)\omega(C_q)=q-1<q$. Lemma~\ref{lem:substitution}
gives the corresponding entries of the table.

The order-$8$ example is $K_2^{\,1}[K_2^{\,2}[K_2^{\,3}]]$, the
order-$12$ example is $K_2^{\,1}[K_2^{\,2}[K_3^{\,3}]]$, and the
order-$16$ example is
$K_2^{\,1}[K_2^{\,2}[K_2^{\,3}[K_2^{\,4}]]]$. All three therefore
avoid $\Pi$ and $\Delta$ and are CIS.

Let $H_{36}$ be the two-colored circulant of order $36$ with lengths
$2,3,6,10,12,14,15,18$ of the first color and all other lengths of
the second color. This circulant is CIS: all maximal cliques and
maximal stable sets of its first chromatic component have size $6$
\cite[Section~3.2]{BGM14}. The vertices $2,0,3,1$, in this order,
induce a $P_4$ in that component and hence a copy of $\Pi$.
Thus a CIS circulant need not be $\Pi$-free.

Figures~\ref{f2}--\ref{f6} show five reduced game trees, their complete
colored graphs, and their circular representations. Each level is
assigned a different player in these examples. The color convention
is $1$ red, $2$ blue, $3$ green, and $4$ magenta. The permutations in
the captions are the mixed-radix reorderings from the proof of
Theorem~\ref{thm:divisibility}; for binary trees they reverse the order
of the binary digits.

\GallaiFigureTwo
\GallaiFigureThree
\GallaiFigureFour
\GallaiFigureFive
\GallaiFigureSix
\FloatBarrier

The level labels need not all be distinct. For example, a binary tree
with three internal levels labeled $1,2,1$ gives
\[
K_2^{\,1}[K_2^{\,2}[K_2^{\,1}]]=(1,3,4;2)_8.
\]
Player $1$ moves at the first and third levels, but never at
successive positions. Equivalently, this coloring is obtained from
$(1,3;2;4)_8$ by merging its first and third colors. The tree remains
reduced, and the graph is still $(\Pi,\Delta)$-free and CIS.

\begin{remark}
The normal forms of positional game structures and their correspondence
with $(\Pi,\Delta)$-free $d$-graphs are studied in \cite{Gur78}; see
also \cite{Gur82,Gur84,Gur09,ABG06}. The uniform trees identify the
structures within this correspondence whose outcome graphs have
circular symmetry.
\end{remark}

The two decompositions play different roles. For Gallai circulants
the factors may be arbitrary two-colored circulants, and the CIS
question remains with those factors. For $(\Pi,\Delta)$-free
circulants the factors are monochromatic, so the branching numbers
and level labels form an intrinsic parameter sequence. They give the
edge colors in the original cyclic labeling through
\eqref{eq:divisibility-colors} and can be recovered by
Corollary~\ref{cor:recognition}.

\FloatBarrier
\section*{Acknowledgements}
This article is an output of a research project (HSE-BR-2025-024)
implemented as part of the Basic Research Program at HSE University.


\begin{thebibliography}{99}

\bibitem{ABG05}
D. V. Andrade, E. Boros, and V. Gurvich,
Even hole free and balanced circulants,
RUTCOR Research Report RRR-08-2005, Rutgers University, 2005.

\bibitem{ABG06}
D. V. Andrade, E. Boros, and V. Gurvich,
On graphs whose maximal cliques and stable sets intersect,
RUTCOR Research Report RRR 17-2006 and DIMACS Technical Report
2006-16, Rutgers University, 2006; published in
B. Goldengorin (ed.), Optimization Problems in Graph Theory,
Springer Optimization and Its Applications 139 (2018), 3--63.

\bibitem{ABG10}
D. V. Andrade, E. Boros, and V. Gurvich,
Not complementary connected and not CIS $d$-graphs form weakly
monotone families,
Discrete Math. 310(5) (2010), 1089--1096.

\bibitem{AFG98}
A. Apartsin, E. Ferapontova, and V. Gurvich,
A circular graph-counterexample to the Duchet kernel conjecture,
Discrete Math. 178 (1998), 229--231.

\bibitem{BBGMP98}
G. Bacs\'o, E. Boros, V. Gurvich, F. Maffray, and M. Preissmann,
On minimal imperfect graphs with circular symmetry,
J. Graph Theory 29(4) (1998), 209--224.

\bibitem{BD98}
S. B\"ocker and A. W. M. Dress,
Recovering symbolically dated, rooted trees from symbolic ultrametrics,
Adv. Math. 138(1) (1998), 105--125.
\url{https://doi.org/10.1006/aima.1998.1743}.

\bibitem{BG93}
E. Boros and V. Gurvich,
When is a circular graph minimal imperfect?,
RUTCOR Research Report RRR 22-1993, Rutgers University, 1993.

\bibitem{BG09}
E. Boros and V. Gurvich,
Vertex- and edge-minimal and locally minimal graphs,
Discrete Math. 309(12) (2009), 3853--3865.

\bibitem{BGM14}
E. Boros, V. Gurvich, and M. Milani\v{c},
On CIS circulants,
Discrete Math. 318 (2014), 78--95.
\url{https://doi.org/10.1016/j.disc.2013.11.015}.

\bibitem{BGM17}
E. Boros, V. Gurvich, and M. Milani\v{c},
On equistable, split, CIS, and related classes of graphs,
Discrete Appl. Math. 216 (2017), 47--66.

\bibitem{BGZ09}
E. Boros, V. Gurvich, and I. Zverovich,
On split and almost CIS-graphs,
Australas. J. Combin. 43 (2009), 163--180.

\bibitem{CE97}
K. Cameron and J. Edmonds,
Lambda composition,
J. Graph Theory 26 (1997), 9--16.

\bibitem{CEL86}
K. Cameron, J. Edmonds, and L. Lov\'asz,
A note on perfect graphs,
Period. Math. Hungar. 17(3) (1986), 173--175.
\url{https://doi.org/10.1007/BF01848646}.

\bibitem{CG83}
F. R. K. Chung and R. L. Graham,
Edge-colored complete graphs with precisely colored subgraphs,
Combinatorica 3 (1983), 315--324.

\bibitem{CLB81}
D. G. Corneil, H. Lerchs, and L. Stewart Burlingham,
Complement reducible graphs,
Discrete Appl. Math. 3 (1981), 163--174.

\bibitem{DHMV15}
E. Dobson, A. Hujdurovi\'c, M. Milani\v{c}, and G. Verret,
Vertex-transitive CIS graphs,
European J. Combin. 44 (2015), 87--98.
\url{https://doi.org/10.1016/j.ejc.2014.09.007}.

\bibitem{DLZ04}
X. Deng, G. Li, and W. Zang,
Proof of Chv\'atal's conjecture on maximal stable sets and maximal
cliques in graphs,
J. Combin. Theory Ser. B 91(2) (2004), 301--325.

\bibitem{DLZ04a}
X. Deng, G. Li, and W. Zang,
Corrigendum to ``Proof of Chv\'atal's conjecture on maximal stable
sets and maximal cliques in graphs'',
J. Combin. Theory Ser. B 94(2) (2005), 352--353.

\bibitem{Eit94}
T. Eiter,
Exact transversal hypergraphs and application to Boolean
$\mu$-functions,
J. Symbolic Comput. 17 (1994), 215--225.

\bibitem{Eit95}
T. Eiter,
Generating Boolean $\mu$-expressions,
Acta Inform. 32 (1995), 171--187.

\bibitem{ER90}
A. Ehrenfeucht and G. Rozenberg,
Theory of 2-structures, part II: Representation through labeled tree
families,
Theoret. Comput. Sci. 70(3) (1990), 305--342.
\url{https://doi.org/10.1016/0304-3975(90)90130-A}.

\bibitem{Ful71}
D. R. Fulkerson,
Blocking and anti-blocking pairs of polyhedra,
Math. Programming 1 (1971), 168--194.

\bibitem{Gal67}
T. Gallai,
Transitiv orientierbare Graphen,
Acta Math. Acad. Sci. Hungar. 18(1--2) (1967), 25--66;
English translation by F. Maffray and M. Preissmann in
J. L. Ram\'irez Alfons\'in and B. A. Reed (eds.),
Perfect Graphs, Wiley, 2002, Chapter 3.

\bibitem{GG11}
M. C. Golumbic and V. Gurvich,
Read-once Boolean functions,
in Y. Crama and P. L. Hammer (eds.),
Boolean Functions: Theory, Algorithms, and Applications,
Cambridge University Press, 2011, Chapter 10, 448--486.

\bibitem{GGSTU93}
A. Gol'berg, V. Gurvich, A. Shapovalov, A. Temkin, and V. Udalov,
Spanning graphs without odd holes and antiholes,
Russian Acad. Sci. Dokl. Math. 47(2) (1993), 278--284.

\bibitem{Gur77}
V. Gurvich,
On repetition-free Boolean functions,
Uspekhi Mat. Nauk 32(1) (1977), 183--184 (in Russian).

\bibitem{Gur78}
V. Gurvich,
Applications of Boolean functions and networks in game theory,
Ph.D. thesis, Moscow Institute of Physics and Technology,
Moscow, USSR, 1978 (in Russian).

\bibitem{Gur82}
V. Gurvich,
On the normal form of positional games,
Soviet Math. Dokl. 25(3) (1982), 572--575.

\bibitem{Gur84}
V. Gurvich,
Some properties and applications of complete edge-chromatic graphs
and hypergraphs,
Soviet Math. Dokl. 30(3) (1984), 803--807.

\bibitem{Gur91}
V. Gurvich,
Criteria for repetition freeness of functions in the algebra of logic,
Russian Acad. Sci. Dokl. Math. 43(3) (1991), 721--726.

\bibitem{Gur09}
V. Gurvich,
Decomposing complete edge-chromatic graphs and hypergraphs. Revisited,
Discrete Appl. Math. 157(14) (2009), 3069--3085.
\url{https://doi.org/10.1016/j.dam.2009.06.026}.

\bibitem{Gur11}
V. Gurvich,
On exact blockers and anti-blockers, $\Delta$-conjecture, and related
problems,
Discrete Appl. Math. 159 (2011), 311--321.

\bibitem{GN26}
V. Gurvich and M. Naumova,
More on discrete convexity,
Discrete Appl. Math. 385 (2026), 146--164.

\bibitem{GS04}
A. Gy\'arf\'as and G. Simonyi,
Edge colorings of complete graphs without tricolored triangles,
J. Graph Theory 46(3) (2004), 211--216.
\url{https://doi.org/10.1002/jgt.20001}.

\bibitem{Hos08}
R. Hoshino,
Independence Polynomials of Circulant Graphs,
Ph.D. thesis, Dalhousie University, 2008.

\bibitem{HW18}
M. Hellmuth and N. Wieseke,
On tree representations of relations and graphs: symbolic ultrametrics
and cograph edge decompositions,
J. Combin. Optim. 36(2) (2018), 591--616.
\url{https://doi.org/10.1007/s10878-017-0111-7}.

\bibitem{Lov72a}
L. Lov\'asz,
Normal hypergraphs and the perfect graph conjecture,
Discrete Math. 2 (1972), 253--267.

\bibitem{Lov72b}
L. Lov\'asz,
A characterization of perfect graphs,
J. Combin. Theory Ser. B 13 (1972), 95--98.

\bibitem{LT26}
Y. Liu and Q. Tang,
The $\Delta$-conjecture for CIS $d$-graphs,
arXiv:2608.08799v2, 2026.
\url{https://arxiv.org/abs/2608.08799}.

\bibitem{McS99}
R. M. McConnell and J. P. Spinrad,
Modular decomposition and transitive orientation,
Discrete Math. 201 (1999), 189--241.

\bibitem{Mor85}
R. H. M\"ohring,
Algorithmic aspects of the substitution decomposition in optimization
over relations, set systems, and Boolean functions,
Ann. Oper. Res. 4 (1985), 195--225.

\bibitem{MS89}
J. H. Muller and J. P. Spinrad,
Incremental modular decomposition,
J. ACM 36(1) (1989), 1--19.

\bibitem{Sei74}
D. Seinsche,
On a property of the class of $n$-colorable graphs,
J. Combin. Theory Ser. B 16 (1974), 191--193.

\bibitem{Sum73}
D. P. Sumner,
Graphs indecomposable with respect to the X-join,
Discrete Math. 6 (1973), 281--298.

\bibitem{TYZ26}
R. Tao, M. Yang, and W. Zang,
How difficult is it to recognize CIS graphs?,
arXiv:2608.11289v1 [cs.DM], 2026.
\url{https://arxiv.org/abs/2608.11289}.

\bibitem{WZZ09}
Y. Wu, W. Zang, and C.-Q. Zhang,
A characterization of almost CIS graphs,
SIAM J. Discrete Math. 23(2) (2009), 749--753.
\url{https://doi.org/10.1137/080723739}.

\bibitem{Zan95}
W. Zang,
Generalizations of Grillet's theorem on maximal stable sets and maximal
cliques in graphs,
Discrete Math. 143(1--3) (1995), 259--268.

\end{thebibliography}
\end{document}